\pdfoutput=1
\RequirePackage{ifpdf}
\ifpdf % We are running pdfTeX in pdf mode
\documentclass[pdftex]{sigma}
\else
\documentclass{sigma}
\fi

\numberwithin{equation}{section}

\newtheorem{Theorem}{Theorem}[section]
\newtheorem{Lemma}[Theorem]{Lemma}
\newtheorem{Assumption}[Theorem]{Assumption}
 { \theoremstyle{definition}
\newtheorem{Definition}[Theorem]{Definition}
\newtheorem{Example}[Theorem]{Example}
 }

\newcommand{\Ga}{\Gamma}
\newcommand{\Up}{\Upsilon}

\newcommand{\ep}{\epsilon}
\newcommand{\la}{\lambda}
\newcommand{\si}{\sigma}

\newcommand{\bC}{\mathbb{C}}
\newcommand{\bE}{\mathbb{E}}

\newcommand{\bL}{\mathbb{L}}
\newcommand{\bP}{\mathbb{P}}
\newcommand{\bQ}{\mathbb{Q}}
\newcommand{\bR}{\mathbb{R}}
\newcommand{\bZ}{\mathbb{Z}}

\newcommand{\cF}{\mathcal{F}}

\newcommand{\cM}{\mathcal{M}}

\newcommand{\cO}{\mathcal{O}}

\newcommand{\cT}{\mathcal{T}}

\newcommand{\Aut}{\mathrm{Aut}}

\newcommand{\ch}{\mathrm{ch}}

\newcommand{\ev}{\mathrm{ev}}
\newcommand{\Ext}{\mathrm{Ext}}
\newcommand{\Hom}{\mathrm{Hom}}

\newcommand{\val}{\mathrm{val}}
\newcommand{\vir}{\mathrm{vir} }

\newcommand{\bh}{ {\mathbf{h}} }

\newcommand{\bw}{ {\mathbf{w}} }
\newcommand{\bx}{ {\mathbf{x}} }

\newcommand{\vGa}{ {\vec{\Ga }} }

\newcommand{\vmu}{ {\vec{\mu }} }

\newcommand{\vd}{\vec{d}}
\newcommand{\vf}{\vec{f}}
\newcommand{\vg}{\vec{g}}

\newcommand{\vs}{\vec{s}}

\newcommand{\Mbar}{\overline{\cM}}

\newcommand{\MgX}{\Mbar_{g,n}(X,\beta)}
\newcommand{\GgX}{G_{g,n}(X,\beta)}

\newcommand{\edge}{ { e\in E(\Ga) \atop (e,v),(e,v')\in F(\Ga)} }
\newcommand{\vone}{ {v\in V^1(\Ga),\, (e,v)\in F(\Ga)} }

\newcommand{\pt}{{\rm point}}

\begin{document}

\allowdisplaybreaks

\newcommand{\arXivNumber}{1407.1370}

\renewcommand{\thefootnote}{}

\renewcommand{\PaperNumber}{048}

\FirstPageHeading

\ShortArticleName{Equivariant Gromov--Witten Invariants of Algebraic GKM Manifolds}

\ArticleName{Equivariant Gromov--Witten Invariants\\ of Algebraic GKM Manifolds\footnote{This paper is a~contribution to the Special Issue on Combinatorics of Moduli Spaces: Integrability, Cohomo\-logy, Quantisation, and Beyond. The full collection is available at \href{http://www.emis.de/journals/SIGMA/moduli-spaces-2016.html}{http://www.emis.de/journals/SIGMA/moduli-spaces-2016.html}}}

\Author{Chiu-Chu Melissa LIU~$^\dag$ and Artan SHESHMANI~$^{\ddag\S}$}

\AuthorNameForHeading{C.-C.M.~Liu and A.~Sheshmani}

\Address{$^\dag$~Department of Mathematics, Columbia University,\\
\hphantom{$^\dag$}~2990 Broadway, New York, NY 10027, USA}
\EmailD{\href{mailto:ccliu@math.columbia.edu}{ccliu@math.columbia.edu}}

\Address{$^\ddag$~Harvard University, Department of Mathematics (CMSA),\\
 \hphantom{$^\ddag$}~20 Garden Street, Cambridge, MA, 02138, USA}
\EmailD{\href{mailto:artan@cmsa.fas.harvard.edu}{artan@cmsa.fas.harvard.edu}}

\Address{$^\S$~Aarhus University, Department of Mathematics,\\
\hphantom{$^\S$}~QGM, Ny Munkegade 118, 8000 Aarhus, Denmark}
\EmailD{\href{mailto:artan@qgm.au.dk}{artan@qgm.au.dk}}

\ArticleDates{Received January 16, 2017, in f\/inal form June 21, 2017; Published online July 01, 2017}

\Abstract{An algebraic GKM manifold is a non-singular algebraic variety equipped with an algebraic action of an algebraic torus, with only f\/initely many torus f\/ixed points and f\/initely many 1-dimensional orbits. In this expository article, we use virtual localization to express equivariant Gromov--Witten invariants of any algebraic GKM manifold (which is not necessarily compact) in terms of Hodge integrals over moduli stacks of stable curves and the GKM graph of the GKM manifold.}


\Classification{14C05; 14D20; 14F05; 14J30; 14N10}

\renewcommand{\thefootnote}{\arabic{footnote}}
\setcounter{footnote}{0}

\section{Introduction}

Gromov--Witten invariants of a projective manifold $X$ are virtual counts of parametrized algebraic curves in $X$, and can be viewed as intersection numbers on moduli spaces of stable maps $\Mbar_{g,n}(X,\beta)$ to $X$. If $X$ is equipped with an algebraic action by a torus $T$, then $T$ acts on the moduli spaces $\Mbar_{g,n}(X,\beta)$, and the Gromov--Witten invariants on $X$ can, by localization, be reduced to the intersection theory on torus f\/ixed substack $\Mbar_{g,n}(X,\beta)^T$ in $\Mbar_{g,n}(X,\beta)$.

An algebraic GKM manifold, named after Goresky--Kottwitz--MacPherson, is a non-singular algebraic variety equipped with an algebraic action of $T$, such that there are f\/initely many 0-dimensional and 1-dimensional orbits. Examples of algebraic GKM manifolds include toric manifolds, Grassmanians, f\/lag manifolds, etc. If $X$ is an algebraic GKM manifold then each connected component of $\Mbar_{g,n}(X,\beta)^T$ is, up to some quasi-f\/inite map, a product of moduli stacks of pointed stable curves, and localization computations reduce Gromov--Witten invariants of~$X$ to Hodge integrals, which are intersection numbers on the moduli stacks of pointed stable curves. This algorithm was f\/irst described by Kontsevich for
genus zero Gromov--Witten invariants of~$\bP^r$ in 1994~\cite{Ko2}, before the construction of virtual fundamental class and the proof of virtual
localization. The moduli spaces $\Mbar_{0,n}(\bP^r,d)$ of genus zero stable maps to $\bP^r$ are proper {\em smooth} DM stacks, so there exists a fundamental class $[\Mbar_{0,n}(\bP^r,d)]\in H_*(\Mbar_{0,n}(\bP^r,d);\bQ)$, and one may apply the classical Atiyah--Bott localization formula~\cite{AtBo} in this case. In \cite[Section~4]{GrPa}, T.~Graber and R.~Pandharipande used their virtual localization formula to derive an explicit formula for all genus Gromov--Witten invariants of~$\bP^r$. (See also K.~Behrend \cite[Section~4]{Be2}.) H.~Spielberg derived a formula of genus zero Gromov--Witten invariants of toric manifolds in his thesis~\cite{Sp}. Localization computations of all genus equivariant Gromov--Witten invariants of toric manifolds can be found in~\cite{Liu}.

The main purpose of this paper is to provide details of the virtual localization calculations of all genus equivariant Gromov--Witten invariants for general algebraic GKM manifolds described on pp.~20--21 of preprint version of~\cite{GrPa}\footnote{Available at \url{http://www.math.ethz.ch/~rahul/loc.ps}.}, for readers who are interested in such details. We do not assume the reader is familiar with GKM manifolds and Gromov--Witten theory: in Section~\ref{sec:GKM}, we def\/ine algebraic GKM manifolds and their GKM graphs, following~\cite{GKM, GZ}; in Section~\ref{sec:GW}, we give a brief review of Gromov--Witten theory. The main computations and formulae are presented in Section~\ref{sec:localization}: we compute all genus equivariant descendant Gromov--Witten invariants of an arbitrary algebraic GKM manifold by virtual localization, and express the answer in terms of Hodge integrals and the GKM graph of the algebraic GKM manifold. Most of Section~\ref{sec:localization} is the straightforward generalization of the $\bP^r$ case discussed in~\cite{Ko2} (genus~0) and
\cite[Section~4]{GrPa}, \cite[Section~4]{Be2} (higher genus); see also \cite[Chapter~27]{HKKPTVVZ}.

\section{Algebraic GKM manifolds} \label{sec:GKM}

In this section, we review the geometry of algebraic GKM manifolds, following \cite{GKM}, and introduce the GKM graph associated to an algebraic GKM manifold, following \cite{GZ}. The GKM graph in this paper can be non-compact since we consider algebraic GKM manifolds which are not necessarily compact. In Section \ref{sec:localization}, we will see that the GKM graph contains all the information needed for computing the Gromov--Witten invariants and the equivariant Gromov--Witten invariants of the GKM manifold.

\subsection{Basic notation}
In this paper, we work over $\bC$.

Let $X$ be a non-singular algebraic variety of dimension $r$. We say that $X$ is an algebraic GKM manifold if it is equipped with an algebraic action of a complex algebraic torus $T=(\bC^*)^m$ with only f\/initely many torus f\/ixed points and f\/initely many 1-dimensional orbits.

Let $N=\Hom(\bC^*, T) \cong \bZ^m$ be the lattice of 1-parameter subgroups of $T$, and let $M=\Hom(T,\bC^*)$ be the lattice of irreducible characters of $T$. Then $M=\Hom(N,\bZ)$ is the dual lattice of~$N$. Let $N_\bR = N\otimes_\bZ \bR$ and $M_\bR = M\otimes_\bZ \bR$, so that they are dual real vector spaces of dimension~$m$. Let $N_\bQ = N\otimes_\bZ \bQ$ and let $M_\bQ = M\otimes_{\bZ}\bQ$. Then $M_\bQ$ can be canonically identif\/ied with $H^2_T(\pt;\bQ)$.

Let
\begin{gather*}
R_T := H^*_T(\{ \pt \};\bQ) =H^*(BT;\bQ)=\bQ[u_1,\ldots,u_m]
\end{gather*}
be the $T$-equivariant cohomology of a point, where $u_i\in H^2_T(BT;\bQ)$. Let $Q_T=\bQ(u_1,\ldots,u_m)$ be the fractional f\/ield of $R_T$.

We make the following assumption on $X$.
\begin{Assumption}\label{assume}\quad
\begin{enumerate}\itemsep=0pt
\item[$1.$] The set $X^T$ of $T$ fixed points in $X$ is non-empty.
\item[$2.$] The closure of a $1$-dimensional orbit is either a complex projective line $\bP^1$ or a complex affine line~$\bC$.
\end{enumerate}
\end{Assumption}
Note that 1) and 2) hold when $X$ is a proper algebraic GKM manifold. Indeed, if $X$ is a~proper algebraic GKM manifold then the closure of any 1-dimensional orbit is $\bP^1$.

\begin{Example}
If $X$ is a non-singular toric variety def\/ined by a f\/inite fan, then $X$ is an algebraic GKM manifold.
\end{Example}

\begin{Example}[the Grassmannian ${\rm Gr}(k,m)$] \label{Gr}
Let ${\rm Gr}(k,m)$ be the set of $k$-dimensional linear subspace of $\bC^m$. It is a nonsingular projective variety of dimension $k(m-k)$. Let $T=(\bC^*)^m$ act on $\bC^m$ by
\begin{gather*}
(t_1,\ldots, t_m)\cdot (z_1,\ldots,z_m)=(t_1 z_1,\ldots, t_m z_m).
\end{gather*}
Given $t\in T$, let $\phi_t\colon \bC^m\to \bC^m$ be def\/ined by $\phi_t(z)=t\cdot z$. Let $T$ act on ${\rm Gr}(k,m)$ by $t\cdot V = \phi_t(V)$, where $V$ is a $k$-dimensional linear subspace of $\bC^m$. Given $J\subset\{1,\ldots,m\}$, let $J^c:= \{1,\ldots,m\}\setminus J$, and def\/ine
\begin{gather*}
\bC^J:=\big\{ (z_1,\ldots, z_m)\in \bC^m\colon z_i=0 \text{ if } i\in J^c\big\} \cong \bC^{|J|}.
\end{gather*}
Note that $\phi_t\big(\bC^J\big)=\bC^J$ for any $t\in T$, $J\subset \{1,\ldots,n\}$.

The torus-f\/ixed points in ${\rm Gr}(k,m)$ are
\begin{gather*}
{\rm Gr}(k,m)^T=\big\{ \bC^J\colon J\subset\{1,\ldots,n\},\, |J|=k\big\}.
\end{gather*}
So there are ${m\choose k}$ torus-f\/ixed points in ${\rm Gr}(k,m)$.

Let $\bC^J$ and $\bC^{J'}$ be distinct $T$-f\/ixed points in ${\rm Gr}(k,m)$. Then $\bC^J\cap\bC^{J'}= \bC^{J\cap J'}$. There is a~torus-f\/ixed line connecting $\bC^J$ and $\bC^{J'}$ if and only of $|J\cap J'|=k-1$. In this case, $|J\cup J'|=k+1$. The $T$-f\/ixed lines in ${\rm Gr}(k,m)$ are
\begin{gather*}
\big\{ \ell_{I,K}\colon I\subset K\subset \{1,\ldots,m\},\, |I|=k-1,\, |K|=k+1\big\},
\end{gather*}
where
\begin{gather*}
\ell_{I,K}=\big\{ V\in {\rm Gr}(k,m)\colon \bC^I\subset V\subset \bC^K\big\} \cong \bP^1.
\end{gather*}
Suppose that $I\subset K\subset \{1,\ldots,m\}$, and $|I|=k-1$, $|K|=k+1$. Then $K=I\cup\{j_1,j_2\}$, where $j_1, j_2\in I^c$. So there are ${m\choose k-1}{m-k+1\choose 2}$ torus-f\/ixed lines in ${\rm Gr}(k,m)$.
\end{Example}

\subsection{GKM graph}\label{GKM-graph}
Let $X$ be an algebraic GKM manifold of dimension $r$, so that $T=(\bC^*)^m$ acts algebraically on~$X$.

Following \cite{GZ}, we def\/ine a graph $\Up$ as follows. Let $V(\Up)$ (resp.~$E(\Up)$) denote the set of vertices (resp.~edges) in~$\Up$.
\begin{enumerate}\itemsep=0pt
\item (Vertices) We assign a vertex $\si$ to each torus f\/ixed point $p_\sigma$ in $X$.
\item (Edges) We assign an edge $\ep$ to each 1-dimensional $O_\ep$ in $X$. Let $\ell_e$ be the closure of $O_\ep$.
\item (Flags) The set of f\/lags in the graph $\Up$ is given by
\begin{gather*}
F(\Up)= \big\{ (\ep,\si)\in E(\Up)\times V(\Up)\colon \si\in \ep\big\}= \big\{(\ep,\si)\in E(\Up)\times V(\Up)\colon p_\si\in \ell_\ep\big\}.
\end{gather*}
\end{enumerate}
The Assumption \ref{assume} can be rephrased in terms of the graph $\Up$.
\begin{Assumption}\quad
\begin{enumerate}\itemsep=0pt
\item[$1.$] $V(\Up)$ is non-empty.
\item[$2.$] Each edge in $E(\Up)$ contains at least one vertex.
\end{enumerate}
\end{Assumption}

Let $E(\Up)_c =\{ \ep\in E(\Up)\colon \ell_\ep \cong \bP^1\}$ be the set of compact edges in~$\Up$. Note that $E(\Up)_c=E(\Up)$ if $X$ is proper.

Given a vertex $\si\in V(\Up)$, we denote by $E_\si$ the set of edges containing $\si$, i.e., $E_\si:=\{\ep \in E(\Up)\colon (\ep,\si)\in F(\Up)\}$. Then $|E_\si|=r$ for all $\si\in V(\Up)$, so $\Up$ is an $r$-valent graph.

Given a f\/lag $(\ep,\si)\in F(\Up)$, let $\bw(\ep,\si) \in M=\Hom(T,\bC^*)$ be the weight of $T$-action on~$T_{p_\si} \ell_\ep$, the tangent line to $\ell_\ep$ at the f\/ixed point $p_\si$, namely,
\begin{gather*}
\bw(\ep,\si) := c_1^T(T_{p_\si}\ell_\ep) \in H^2_T(p_\si;\bZ)\cong M.
\end{gather*}
This gives rise to a map $\bw\colon F(\Up)\to M$ satisfying the following
properties.
\begin{enumerate}\itemsep=0pt
\item (GKM hypothesis) Given any $\si\in V(\Up)$, and any two distinct edges $\ep, \ep'\in E_\si$, $\bw(\ep,\si)$ and $\bw(\ep',\si)$ are linearly independent in $M_\bR$.

\item Any edge $\ep\in E_\si$ connecting the vertices $\si,\si'\in V(\Up)$ satisf\/ies the property that:
\begin{enumerate}\itemsep=0pt
\item $\bw(\ep,\si)+\bw(\ep,\si')=0$.
\item Let $E_\si= \{ \ep_1,\ldots, \ep_r\}$, where $\ep_r=\ep$. For any $\ep_i\in E_\si$ there exists $\ep_i'\in E_{\si'}$ and $a_i\in \bZ$ such that
\begin{gather*}
\bw(\ep_i',\si') =\bw(\ep_i,\si)- a_i \bw(\ep,\si).
\end{gather*}
In particular, $\ep'_r=\ep_r=\ep$ and $a_r=2$.
\end{enumerate}
\end{enumerate}
Let $\ep$ be as in item 2 above. The normal bundle of $\ell_\ep \cong \bP^1$ in $X$ is given by
\begin{gather*}
N_{\ell_\ep/X}\cong L_1\oplus\cdots \oplus L_{r-1},
\end{gather*}
where $L_i$ is a degree $a_i$ $T$-equivariant line bundle over $\ell_\ep$ such that the weights of the $T$-action on the f\/ibers $(L_i)_{p_\si}$ and $(L_i)_{p_{\si'} }$ are $\bw(\ep_i,\si)$ and $\bw(\ep_i',\si')$, respectively. The map $\bw\colon F(\Up)\to M$ is called the {\em axial function}.

\begin{Example}[${\rm Gr}(k,m)$] The GKM graph of ${\rm Gr}(k,m)$ is a $k(m-k)$-valent graph $\Up$ such that
\begin{gather*}
V(\Up) = \big\{\si_J\colon J\subset \{1,\ldots,n\},|J|=k\big\},\\
E(\Up) = E(\Up)_c =\big\{ \ep_{I,K}\colon I\subset K\subset\{1,\ldots,m\}, |I|=k-1, |K|=k+1\big\},\\
F(\Up)= \big\{ (\ep_{I,K},\si_J)\in E(\Up)\times V(\Up)\colon I\subset J\subset K\big\},\\
\bw\big(\ep_{IK},\si_{I\cup \{j_{1}\}}\big) = -\bw\big(\ep_{IK},\si_{I\cup \{j_{2}\}}\big) = u_{j_{2}}-u_{j_{1}},\qquad j_{1}, j_{2}\in I^c, \qquad K=I\cup \{j_{1},j_{2}\}.
\end{gather*}
\end{Example}

We def\/ine the 1-skeleton of $X$ to be the union of 1-dimensional orbit closures:
\begin{gather*}
X^1:=\bigcup_{\ep\in E(\Up)} \ell_\ep.
\end{gather*}
The formal completion $\hat{X}$ of $X$ along the 1-skeleton $X^1$ (def\/ined on p.~194 of~\cite{Ha}), together with the $T$-action, can be reconstructed from the graph $\Up$ and $\bw\colon F(\Up)\to M$. We call $(\Up,\bw)$ the GKM graph of $X$ with the $T$-action. If $\rho\colon T'\to T$ is a homomorphism between complex algebraic tori, then~$T'$ acts on~$X$ by $t'\cdot x = \rho(t')\cdot x$, where $t'\in T'$, $\rho(t')\in T$, $x\in X$. If the 0-dimensional and 1-dimensional orbits of this $T'$-action coincide with those of the $T$-action, then the GKM graph of $X$ with this $T'$-action is given by $(\Up,\rho^*\circ \bw)$, where $\rho^*\colon M= \Hom(T,\bC^*)\to \Hom(T',\bC^*)$.

\section{Gromov--Witten theory} \label{sec:GW}

In this section, we give a brief review of the Gromov--Witten theory and the equivariant Gromov--Witten theory.

\subsection{Moduli space of stable curves and Hodge integrals}
\label{sec:hodge}

An $n$-pointed, genus $g$ prestable curve
is a connected algebraic curve $C$ of
arithmetic genus $g$ together with
$n$ ordered marked points $x_1,\ldots, x_n\in C$,
where $C$ has at most nodal singularities,
and $x_1,\ldots,x_n$ are distinct smooth points.
An $n$-pointed, genus $g$ prestable
curve $(C,x_1,\ldots,x_n)$ is {\em stable} if its automorphism
group is f\/inite, or equivalently,
\begin{gather*}
\Hom_{\cO_C}\big(\Omega_C(x_1+\cdots+x_n), \cO_C\big) =0.
\end{gather*}

Let $\Mbar_{g,n}$ be the moduli space of $n$-pointed, genus $g$ stable curves, where $n, g$ are non-negative integers. We assume that $2g-2+n>0$, so that $\Mbar_{g,n}$ is nonempty. Then $\Mbar_{g,n}$ is a proper smooth Deligne--Mumford stack of dimension $3g-3+n$ \cite{DeMu, KnMu, Kn2, Kn3}. The tangent space of~$\Mbar_{g,n}$ at a~moduli point $[(C,x_1,\ldots,x_n)]\in \Mbar_{g,n}$ is given by
\begin{gather*}
\Ext^1_{\cO_C}\big(\Omega_C(x_1+\cdots+x_n), \cO_C\big).
\end{gather*}
Since $\Mbar_{g,n}$ is a proper smooth Deligne--Mumford stack, there is a fundamental class
\begin{gather*}
\big[\,\Mbar_{g,n}\big] \in A_{3g-3+n}\big(\Mbar_{g,n};\bQ\big),
\end{gather*}
which allows us to def\/ine
\begin{gather*}
\int_{\Mbar_{g,n}}\colon \ A^*\big(\Mbar_{g,n};\bQ\big)\longrightarrow \bQ, \qquad \alpha \mapsto \deg\big(\alpha\cap \big[\Mbar_{g,n}\big]\big).
\end{gather*}

We now introduce some classes in $A^*(\Mbar_{g,n};\bQ)$. There is a forgetful morphism $\pi\colon \Mbar_{g,n+1}\to \Mbar_{g,n}$ given by forgetting the $(n+1)$-th marked point (and contracting the unstable irreducible component if there is one):
\begin{gather*}
[(C,x_1,\ldots,x_n,x_{n+1})]\mapsto \big[\big(C^{\rm st}, x_1,\ldots, x_n\big)\big],
\end{gather*}
where $(C^{\rm st}, x_1,\ldots,x_n)$ is the stabilization of the prestable curve $(C,x_1,\ldots,x_n)$. $\pi\colon \Mbar_{g,n+1}\to \Mbar_{g,n}$ can be identif\/ied with the universal curve over $\Mbar_{g,n}$.
\begin{itemize}\itemsep=0pt
\item ($\lambda$ classes) Let $\omega_\pi$ be the relative dualizing sheaf of $\pi\colon \Mbar_{g,n+1}\to \Mbar_{g,n}$. The Hodge bundle $\bE=\pi_*\omega_\pi$ is a rank $g$ vector bundle over $\Mbar_{g,n}$ whose f\/iber over the moduli point $[(C,x_1,\ldots,x_n)]\in \Mbar_{g,n}$
is $H^0(C,\omega_C)$, the space of sections of the dualizing sheaf $\omega_C$ of the curve $C$. The $\la$ classes are def\/ined by
\begin{gather*}
\la_j=c_j(\bE)\in A^j\big(\Mbar_{g,n};\bQ\big).
\end{gather*}
\item ($\psi$ classes) The $i$-th marked point $x_i$ gives rise a section $s_i\colon \Mbar_{g,n}\to \Mbar_{g,n+1}$ of the universal curve. Let $\bL_i=s_i^*\omega_\pi$ be the line bundle over $\Mbar_{g,n}$ whose f\/iber over the moduli point $[(C,x_1,\ldots,x_n)]\in \Mbar_{g,n}$ is the
cotangent line $T_{x_i}^*C$ of $C$ at $x_i$. The $\psi$ classes are def\/ined by
\begin{gather*}
\psi_i=c_1(\bL_i)\in A^1\big(\Mbar_{g,n};\bQ\big).
\end{gather*}
\end{itemize}

{\em Hodge integrals} are top intersection numbers of $\lambda$ classes and $\psi$ classes:
\begin{gather}\label{eqn:hodge}
\int_{\Mbar_{g,n}}\psi_1^{a_1}\cdots \psi_n^{a_n} \la_1^{k_1} \cdots \la_g^{k_g} \in \bQ.
\end{gather}
By def\/inition, \eqref{eqn:hodge} is zero unless
\begin{gather*}
a_1+\cdots + a_n + k_1 + 2k_2 + \cdots + g k_g = 3g-3+n.
\end{gather*}

Using Mumford's Grothendieck--Riemann--Roch calculations in \cite{Mu}, Faber proved, in \cite{Fa}, that general Hodge integrals can be uniquely reconstructed from the $\psi$ integrals (also known as {\em descendant integrals}):
\begin{gather}\label{eqn:descendant}
\int_{\Mbar_{g,n}}\psi_1^{a_1}\cdots \psi_n^{a_n}.
\end{gather}
The descendant integrals can be computed recursively by Witten's conjecture which asserts that the $\psi$ integrals \eqref{eqn:descendant} satisfy a system of dif\/ferential equations known as the KdV equations \cite{Wi}. The KdV equations and the string equation determine all the~$\psi$ integrals~\eqref{eqn:descendant} from the initial value $\int_{\Mbar_{0,3}} 1=1$. For example, from the initial value $\int_{\Mbar_{0,3}}1=1$ and the string equation, one can derive the following formula of genus $0$ descendant integrals:
 \begin{gather*}%\label{eqn:genus-zero-psi}
\int_{\Mbar_{0,n}}\psi_1^{a_1}\cdots\psi_n^{a_n}=\frac{(n-3)!}{a_1!\cdots a_n!},
\end{gather*}
where $a_1+\cdots+a_n=n-3$ \cite[Section~3.3.2]{Ko2}.

The Witten's conjecture was f\/irst proved by Kontsevich in \cite{Ko1}. By now, Witten's conjecture has been reproved many times (Okounkov--Pandharipande \cite{OP1}, Mirzakhani \cite{Mi}, Kim--Liu \cite{KiL}, Kazarian--Lando \cite{KaL}, Chen--Li--Liu \cite{CLL}, Kazarian \cite{Ka}, Mulase--Zhang \cite{MuZ}, etc.).

\subsection{Moduli of stable maps} \label{sec:stable-maps}

Let $X$ be a nonsingular projective or quasi-projective variety, and let $\beta \in H_2(X;\bZ)$. An $n$-pointed, genus $g$, degree $\beta$ prestable map to $X$ is a morphism $f\colon (C,x_1,\ldots,x_n)\to X$, where $(C,x_1,\ldots,x_n)$ is an $n$-pointed, genus~$g$ prestable curve, and $f_*[C]= \beta$. Two prestable maps
\begin{gather*}
f\colon \ (C,x_1,\ldots,x_n)\to X,\qquad f'\colon \ (C', x'_1,\ldots, x'_n)\to X
\end{gather*}
are isomorphic if there exists an isomorphism $\phi\colon (C,x_1,\ldots,x_n)\to (C',x'_1,\ldots,x'_n)$ of $n$-pointed prestable curves such that $f=f'\circ \phi$. A prestable map $f\colon (C,x_1,\ldots,x_n)\to X$ is {\em stable} if its automorphism group is f\/inite. The notion of stable maps was introduced by Kontsevich~\cite{Ko2}.

The moduli space $\Mbar_{g,n}(X,\beta)$ of $n$-pointed, genus $g$, degree $\beta$ stable maps to $X$ is a Deligne--Mumford stack which is proper when $X$ is projective \cite{BeMa}.

\subsection{Obstruction theory and virtual fundamental classes} \label{sec:GWdeform}

The tangent space $T^1$ and the obstruction space $T^2$ at a moduli point $[f\colon (C, x_1,\ldots, x_n)\to X] \in \Mbar_{g,n}(X,\beta)$ f\/it in the {\em tangent-obstruction exact sequence}:
\begin{gather*}
0 \to \Ext^0_{\cO_C}\big(\Omega_C(x_1+\cdots + x_n) , \cO_C\big)\to H^0(C,f^*TX) \to T^1 \nonumber\\
\hphantom{0} {}\to \Ext^1_{\cO_C}\big(\Omega_C(x_1+\cdots + x_n), \cO_C\big) \to H^1(C,f^*TX)\to T^2 \to 0,%\label{eqn:tangent-obstruction}
\end{gather*}
where
\begin{itemize}\itemsep=0pt
\item $\Ext^0_{\cO_C}(\Omega_C(x_1+\cdots+x_n),\cO_C)$ is the space of inf\/initesimal automorphisms of the domain $(C, x_1,\ldots, x_n)$,
\item $\Ext^1_{\cO_C}(\Omega_C(x_1+\cdots + x_n), \cO_C)$ is the space of inf\/initesimal
deformations of the domain $(C, x_1, \ldots, x_n)$,
\item $H^0(C,f^*TX)$ is the space of inf\/initesimal deformations of the map $f$, and
\item $H^1(C,f^*TX)$ is the space of obstructions to deforming the map $f$.
\end{itemize}
$T^1$ and $T^2$ form sheaves $\cT^1$ and $\cT^2$ on the moduli space $\MgX$.

We say $X$ is {\em convex} if $H^1(C,f^*TX)=0$ for all genus $0$ stable maps $f$. Projective spaces $\bP^n$, or more generally, generalized f\/lag varieties $G/P$, are examples of convex varieties. When $X$ is convex and $g=0$, the obstruction sheaf $\cT^2=0$, and the moduli space $\Mbar_{0,n}(X,\beta)$ is a {\em smooth} Deligne--Mumford stack.

In general, $\MgX$ is a {\em singular} Deligne--Mumford stack equipped with a {\em perfect obstruction theory}: there is a two term complex of locally free sheaves $E \to F$ on $\Mbar_{g,n}(X,\beta)$ such that
\begin{gather*}
0\to \cT^1 \to F^\vee \to E^\vee \to \cT^2 \to 0
\end{gather*}
is an exact sequence of sheaves (see \cite{BeFa} for the complete def\/inition of a perfect obstruction theory). The {\em virtual dimension} $d^\vir$ of $\Mbar_{g,n}(X,\beta)$ is the rank of the virtual tangent bundle $T^\vir = F^\vee - E^\vee$,
\begin{gather*}%\label{eqn:virtual-dim}
d^\vir = \int_\beta c_1(TX) +(\dim X-3)(1-g) +n.
\end{gather*}

Suppose that $\MgX$ is {\em proper}. Then there is a {\em virtual fundamental class}
\begin{gather*}
\big[\,\MgX\big]^\vir\in A_{d^\vir}\big(\MgX;\bQ\big).
\end{gather*}
The virtual fundamental class has been constructed
by Li--Tian \cite{LiTi1}, Behrend--Fantechi \cite{BeFa}, Behrend \cite{Be1} in algebraic Gromov--Witten theory. The virtual fundamental class allows us to def\/ine
\begin{gather*}
\int_{[\MgX]^\vir}\colon \ A^*\big(\MgX;\bQ\big) \longrightarrow \bQ,\qquad \alpha \mapsto \deg\big(\alpha\cap\big[\,\MgX\big]^\vir\big).
\end{gather*}

\subsection{Gromov--Witten invariants}\label{sec:GWinvariants}
Let $X$ be a nonsingular projective variety. Gromov--Witten invariants are rational numbers def\/ined by applying
\begin{gather*}
\int_{[\MgX]^\vir}\colon \ A^*\big(\MgX\big)\to \bQ
\end{gather*}
to certain classes in $A^*(\MgX)$.

Let $\ev_i\colon \Mbar_{g,n}(X,\beta)\to X$ be the evaluation at the $i$-th marked point: $\ev_i$ sends $[f\colon (C,x_1,\ldots$, $x_n)\to X]\in \Mbar_{g,n}(X,\beta)$ to $f(x_i)\in X$. Given $\gamma_1,\ldots,\gamma_n \in A^*(X)$, def\/ine
\begin{gather}\label{eqn:primaryGW}
\langle \gamma_1,\ldots, \gamma_n\rangle^X_{g,\beta} = \int_{ [\MgX]^\vir } \ev_1^*\gamma_1 \cup \cdots \cup \ev_n^*\gamma_n \in \bQ.
\end{gather}
These are known as the {\em primary} Gromov--Witten invariants of $X$. More generally, we may also view $[\MgX]^\vir$ as a class in~$H_{2d^{\vir}}(\MgX)$. Then \eqref{eqn:primaryGW} is def\/ined for ordinary cohomology classes $\gamma_1,\ldots,\gamma_n\in H^*(X)$,
including odd cohomology classes which do not come from $A^*(\MgX)$.

Let $\pi\colon \Mbar_{g,n+1}(X,\beta)\!\to\! \Mbar_{g,n}(X,\beta)$ be the universal curve. For $i=1,{\ldots},n$, let $s_i\colon \Mbar_{g,n}(X,\beta)\!$ $\to \Mbar_{g,n+1}(X,\beta)$, be the section which corresponds to the $i$-th marked point. Let $\omega_\pi\to \Mbar_{g,n+1}(X,\beta)$ be the relative dualizing sheaf of $\pi$, and let $\bL_i=s_i^*\omega_\pi$ be the line bundle over $\Mbar_{g,n}(X,\beta)$ whose f\/iber at the moduli point $[f\colon (C,x_1,\ldots,x_n)\to X]\in \Mbar_{g,n}(X,\beta)$ is the cotangent line $T^*_{x_i}C$ at the $i$-th marked point~$x_i$. The $\psi$-classes are def\/ined to be
\begin{gather*}
\psi_i:=c_1(\bL_i)\in A^1\big(\MgX\big),\qquad i=1,\ldots,n.
\end{gather*}
We use the same notation $\psi_i$ to denote the corresponding classes in the ordinary cohomology group $H^2(\MgX)$.

Genus $g$, degree $\beta$ {\em descendant} Gromov--Witten invariants of $X$ are def\/ined by
\begin{gather}\label{eqn:descendantGW}
\langle \tau_{a_1}(\gamma_1)\cdots \tau_{a_n}(\gamma_n)\rangle_{g,\beta}^X
:= \int_{ [\MgX]^\vir } \ev_1^*\gamma_1\cup \psi_1^{a_1} \cup \cdots
\cup \ev_n^*\gamma_n \cup \psi_n^{a_n}
\in \bQ.
\end{gather}
Suppose that $\gamma_i\in H^{d_i}(X)$. Then \eqref{eqn:descendantGW} is zero unless
\begin{gather*}%\label{eqn:nonzero-condition}
\sum_{i=1}^n (d_i + 2a_i-2) =2\bigg(\int_\beta c_1(TX) + (\dim X-3)(1-g)\bigg).
\end{gather*}

\subsection{Equivariant Gromov--Witten invariants} \label{sec:equivariant-GW}
Let $X$ be a non-singular projective or quasi-projective algebraic variety, equipped with an algebraic action of $T=(\bC^*)^m$. Then~$T$ acts on~$\MgX$ by
\begin{gather*}
t\cdot [f\colon (C, x_1,\ldots,x_n)\to X]\mapsto [t\cdot f\colon (C,x_1,\ldots,x_n)\to X],
\end{gather*}
where $(t\cdot f)(z)= t\cdot f(z)$, $z\in C$. The evaluation maps
$\ev_i\colon \MgX\to X$ are $T$-equivariant and induce $\ev_i^*\colon A^*_T(X;\bQ)\to A^*_T(\MgX;\bQ)$.

\subsubsection[Def\/intion when $\MgX$ is proper]{Def\/intion when $\boldsymbol{\MgX}$ is proper}

Suppose that $\MgX$ is proper, so that there are virtual fundamental classes
\begin{gather*}
\big[\,\MgX\big]^\vir \in A_{d^\vir}\big(\MgX;\bQ\big),\qquad \big[\,\MgX\big]^{\vir}_T \in A_{d^\vir}^T\big(\MgX;\bQ\big),
\end{gather*}
where
\begin{gather*}
d^\vir= \int_\beta c_1(TX)+(\dim X-3)(1-g)+n.
\end{gather*}
Given $\gamma_i\in A^{d_i}(X;\bQ)=H^{2d_i}(X;\bQ)$ and $a_i\in \bZ_{\geq 0}$, def\/ine $\langle \tau_{a_i}(\gamma_1)\cdots \tau_{a_n}(\gamma_n)\rangle^X_{g,\beta}$ as in Section~\ref{sec:GWinvariants}:
\begin{gather}\label{eqn:nonequivariantGW}
\langle\tau_{a_1}(\gamma_1)\cdots \tau_{a_n}(\gamma_n)\rangle^X_{g,\beta}
=\int_{[\MgX]^\vir}\prod_{i=1}^n \big(\ev_i^*\gamma_i\cup \psi_i^{a_i}\big) \in \bQ.
\end{gather}
By def\/inition, \eqref{eqn:nonequivariantGW} is zero unless $\sum\limits_{i=1}^n (d_i+a_i) = d^\vir$. In this case,
\begin{gather*}%\label{eqn:alphaT}
\langle \tau_{a_1}(\gamma_1)\cdots \tau_{a_n}(\gamma_n)\rangle^X_{g,\beta}
=\int_{[\MgX]^{\vir}_T}\prod_{i=1}^n \big(\ev_i^*\gamma_i^T\cup \big(\psi_i^T\big)^{a_i}\big),
\end{gather*}
where $\gamma_i^T\in A^{d_i}_T(X)$ is any $T$-equivariant lift of $\gamma_i\in A^{d_i}(X)$, and $\psi_i^T\in A^1_T(\MgX)$ is any $T$-equivariant lift of $\psi_i\in A^1(\MgX)$.

In this paper, we f\/ix a choice of $\psi_i^T$ as follows. A stable map $f\colon (C,x_1,\ldots,x_n)\to X$ induces $\bC$-linear maps $T_{x_i}C\to T_{f(x_i)}X$ for $i=1,\ldots,n$. This gives rise to $\bL_i^\vee \to \ev_i^*TX$. The $T$-action on $X$ induces a $T$-action on $TX$, so that $TX$ is a $T$-equivariant vector bundle over $X$, and $\ev_i^*TX$ is a $T$-equivariant vector bundle over~$\MgX$.

We def\/ine
\begin{gather*}
\psi_i^T=c_1^T(\bL_i)\in A^1_T\big(\MgX\big), \qquad i=1,\ldots,n.
\end{gather*}
Then $\psi_i^T$ is a $T$-equivariant lift of $\psi_i=c_1(\bL_i)\in A^1(\MgX)$.

Given $\gamma_i^T\in A_T^{d_i}(X;\bQ)$,
we def\/ine genus $g$, degree $\beta$ $T$-equivariant descendant Gromov--Witten invariants
\begin{gather*}%\label{eqn:equivariantGW}
\big\langle\tau_{a_1}\big(\gamma_1^T\big),\dots,\tau_{a_n}\big(\gamma_n^T\big)\big\rangle_{g,\beta}^{X_T}
 :=\int_{[\MgX]^{\vir}_T} \prod_{i=1}^n \big(\ev_i^*\gamma_i^T \cup\big(\psi_i^T\big)^{a_i}\big)\\
\hphantom{\big\langle\tau_{a_1}\big(\gamma_1^T\big),\dots,\tau_{a_n}\big(\gamma_n^T\big)\big\rangle_{g,\beta}^{X_T}
 :=}{} \in \bQ[u_1,\ldots,u_m]\left(\sum_{i=1}^n (d_i+a_i)-d^\vir\right),
\end{gather*}
where $\bQ[u_1,\ldots,u_m](k)$ is the space of degree $k$ homogeneous polynomials in $u_1,\ldots,u_m$ with rational coef\/f\/icients. In particular,
\begin{gather*}
\big\langle\tau_{a_1}\big(\gamma_1^T\big),\dots,\tau_{a_n}\big(\gamma_n^T\big)\big\rangle_{g,\beta}^{X_T}
=\begin{cases}
0, & \displaystyle \sum\limits_{i=1}^n (d_i+a_i) <d^\vir,\\
\langle\tau_{a_1}(\gamma_1),\dots,\tau_{a_n}(\gamma_n)\rangle_{g,\beta}^X\in \bQ, &
\displaystyle \sum\limits_{i=1}^n (d_i+a_i) = d^\vir,
\end{cases}
\end{gather*}
where $\gamma_i\in A^{d_i}(X;\bQ)$ is the image of $\gamma_i^T$ under $A^{d_i}_T(X;\bQ)\to A^{d_i}(X;\bQ)$.

{\sloppy Let $\MgX^T\subset \MgX$ be the substack of $T$-f\/ixed points, and let $i\colon \MgX^T\to \MgX$ be the inclusion. Let $N^\vir$ be the virtual normal bundle of substack $\MgX^T$ in~$\MgX$; in general, $N^\vir$ has dif\/ferent ranks on dif\/ferent connected components of $\MgX^T$. By virtual localization,
\begin{gather}\label{eqn:virtual-localization}
\int_{[\MgX]^\vir_T} \prod_{i=1}^n\big(\ev_i^*\gamma_i^T\cup \big(\psi_i^T\big)^{a_i}\big) =\int_{[\MgX^T]^\vir_T}\frac{i^*\prod\limits_{i=1}^n\big(\ev_i^*\gamma_i^T\cup \big(\psi_i^T\big)^{a_i}\big)}{e^T(N^\vir)}.
\end{gather}

}

\subsubsection[Def\/inition when $\MgX$ is not proper]{Def\/inition when $\boldsymbol{\MgX}$ is not proper}
Suppose that $\MgX$ is not proper, but $\MgX^T$ is. Then the left hand side of \eqref{eqn:virtual-localization} is not def\/ined, but the right hand side of \eqref{eqn:virtual-localization} is. In this case, we use the right hand side of \eqref{eqn:virtual-localization} to {\em define} $T$-equivariant Gromov--Witten invariants:
\begin{gather}\label{eqn:residue}
\big\langle \tau_{a_1}\big(\gamma_1^T\big),\ldots,\tau_{a_n}\big(\gamma_n^T\big)\big\rangle^{X_T}_{g,\beta}
:=\int_{[\MgX^T]^{\vir}_T}\frac{i^*\prod\limits_{i=1}^n\big(\ev_i^*\gamma_i^T\cup \big(\psi_i^T\big)^{a_i}\big)}{e^T(N^\vir)}\in Q_T.
\end{gather}
When $\MgX$ is not proper, the right hand side of \eqref{eqn:residue} is a rational function (instead of a polynomial) in $u_1,\ldots,u_m$. It can be nonzero when $\sum\limits_{i=1}^n (d_i+a_i)< d^\vir$, and does not have a~nonequivariant limit (obtained by setting $u_i=0$) in general.

\section{Virtual Localization} \label{sec:localization}

In this section, we compute all genus equivariant descendant Gromov--Witten invariants of any algebraic GKM manifold by virtual localization. This generalizes the toric case in \cite[Section~5]{Liu}.

Let $X$ be an algebraic GKM manifold of dimension $r$, with an algebraic action of $T=(\bC^*)^m$, and let $\Up$ be the corresponding GKM graph.

\subsection{Torus f\/ixed points and graph notation}\label{sec:graph-notation}
In this subsection, we describe the $T$-f\/ixed points in~$\MgX$. Following Kontsevich~\cite{Ko2}, given a stable map $f\colon (C,x_1,\ldots,x_n)\to X$ such that
\begin{gather*}
[f\colon (C,x_1,\ldots,x_n)\to X] \in \MgX^T,
\end{gather*}
we will associate a decorated graph $\vGa$.

We f\/irst give a formal def\/inition.
\begin{Definition}\label{df:GgX}
A decorated graph $\vGa=\big(\Ga, \vf, \vd, \vg, \vs\big)$ for $n$-pointed, genus $g$, degree $\beta$ stable maps to $X$ consists of the following data.
\begin{enumerate}\itemsep=0pt
\item $\Ga$ is a compact, connected 1-dimensional CW complex. We denote the set of vertices (resp.~edges) in $\Ga$ by $V(\Ga)$ (resp.~$E(\Ga)$). Let
\begin{gather*}
F(\Gamma)=\{ (e,v)\in E(\Ga)\times V(\Ga)\,|\, v\in e\}
\end{gather*}
be the set of f\/lags in $\Gamma$.
\item The {\em label map} $\vf\colon V(\Ga)\cup E(\Ga)\to V(\Up)\cup E(\Up)_c$ sends a vertex $v\in V(\Ga)$ to a~vertex $\si_v \in V(\Up)$, and
sends an edge $e\in E(\Ga)$ to an edge $\ep_e \in E(\Up)_c$. Moreover, $\vf$ def\/ines a~map from the graph $\Ga$ to the graph $\Up$: if $(e,v)$ is a f\/lag in $\Gamma$ then $(\ep_e, \si_v)$ is a f\/lag in $\Up$.

\item The {\em degree map} $\vd\colon E(\Ga)\to \bZ_{>0}$ sends an edge $e\in E(\Ga)$ to a positive integer~$d_e$.

\item The {\em genus map} $\vg\colon V(\Ga)\to \bZ_{\geq 0}$ sends a vertex $v\in V(\Ga)$ to a non-negative integer $g_v$.

\item The {\em marking map} $\vs\colon \{1,2,\ldots,n\}\to V(\Ga)$ is def\/ined if $n>0$.
\end{enumerate}

The above maps satisfy the following two constraints:
\begin{enumerate}\itemsep=0pt
\item[(i)] (topology of the domain) $\sum\limits_{v\in V(\Ga)} g_v + |E(\Ga)| - |V(\Ga)| +1 = g$.
\item[(ii)] (topology of the map) $\sum\limits_{e\in E(\Ga)} d_e[\ell_{\ep_e}] =\beta$.
\end{enumerate}

Let $\GgX$ be the set of all decorated graphs $\vGa=\big(\Ga,\vf, \vd,\vg,\vs\big)$ satisfying the above constraints.
\end{Definition}

We now describe the geometry and combinatorics of a stable map $f\colon (C,x_1,\ldots,x_n)\to X$ which represents a $T$-f\/ixed point in~$\MgX$.

For any $t\in T$, there exists an automorphism $\phi_t$ of $(C,x_1,\ldots,x_n)$ such that $t\cdot f(z)= f\circ\phi_t(z)$ for any $z\in C$. Let $C'$ be an irreducible component of $C$, and let $f'=f|_{C'}\colon C'\to X$. There are two possibilities:
\begin{enumerate}\itemsep=0pt
\item[] \hspace*{-10mm} Case~1: $f'$ is a constant map, and $f(C')=\{ p_\si\}$, where $p_\si$ is a f\/ixed point in $X$ associated to some $\si\in V(\Up)$.
\item[] \hspace*{-9.5mm} Case~2: $C'\cong \bP^1$ and $f(C') =\ell_\ep$, where $\ell_\ep$ is a $T$-invariant $\bP^1$ in $X$ associated to some $\ep\in E(\Up)_c$.
\end{enumerate}

We def\/ine a decorated graph $\vGa$ associated to $f\colon (C,x_1,\ldots,x_n)\to X$ as follows.
\begin{enumerate}\itemsep=0pt
\item (Vertices) We assign a vertex $v$ to each connected component $C_v$ of $f^{-1}(X^T)$.
\begin{enumerate}\itemsep=0pt
\item (label) $f(C_v)=\{p_\si\}$ for some $\si\in V(\Up)$;
we def\/ine $\vf(v)=\si_v=\si$.

\item (genus) $C_v$ is a curve or a point. If $C_v$ is a curve
then we def\/ine $\vg(v)=g_v$ to be the arithmetic
genus of $C_v$; if $C_v$ is a point then we def\/ine
$\vg(v)=g_v=0$.
\item (marking) For $i=1,\ldots,n$, def\/ine
$\vs(i)=v$ if $x_i\in C_v$.
\end{enumerate}

\item (Edges)
For any $\ep \in E(\Up)$, let
$O_\ep \cong \bC^*$ be the 1-dimensional
orbit whose closure is $\ell_\ep$.
Then
\begin{gather*}
X^1\setminus X^T = \bigsqcup_{\ep\in E(\Up)} O_\ep,
\end{gather*}
where the right hand side is a disjoint union of connected components. We assign an edge~$e$ to each connected component $O_e\cong \bC^*$ of $f^{-1}\big(X^1\setminus X^T\big)$.
\begin{enumerate}\itemsep=0pt
\item (label) Let $C_e\cong\bP^1$ be the closure of $O_e$. Then $f(C_e)=\ell_\ep$ for some $\ep\in E(\Up)_c$; we def\/ine $\vf(e)=\ep_e=\ep$.
\item (degree) We def\/ine $\vd(e)=d_e$ to be the degree of the map $f|_{C_e}\colon C_e\cong \bP^1\to \ell_\ep\cong \bP^1$.
\end{enumerate}

\item (Flags) The set of f\/lags in the graph $\Ga$ is def\/ined by
\begin{gather*}
F(\Ga)=\{(e,v)\in E(\Ga)\times V(\Ga)\,|\, C_e\cap C_v\neq \varnothing\}.
\end{gather*}
\end{enumerate}
The above 1), 2), 3) def\/ine a decorated graph $\vGa=\big(\Ga, \vf,\vd,\vg,\vs\big)$ satisfying the constraints~(i) and~(ii) in Def\/inition~\ref{df:GgX}.
Therefore $\vGa\in \GgX$. This gives a map from $\MgX^T$ to the discrete set~$\GgX$. Let $\cF_\vGa\subset \MgX^T$ denote the preimage of~$\vGa$.
Then
\begin{gather*}
\MgX^T=\bigsqcup_{\vGa\in \GgX}\cF_\vGa,
\end{gather*}
where the right hand side is a disjoint union of connected components. We next describe the f\/ixed locus $\cF_\vGa$ associated to each decorated graph $\vGa\in \GgX$. For later convenience, we introduce some def\/initions.

\begin{Definition}\label{df:unstable}
Given a vertex $v\in V(\Ga)$, we def\/ine
\begin{gather*}
E_v=\{e\in E(\Ga)\,|\, (e,v)\in F(\Ga)\}
\end{gather*}
the set of edges emanating from $v$, and def\/ine $S_v=\vs^{-1}(v)\subset \{1,\ldots,n\}$. The valency of $v$ is given by $\val(v)= |E_v|$. Let $n_v=|S_v|$ be the number of marked points contained in~$C_v$. We say a vertex is {\em stable} if $2g_v-2 + \val(v) + n_v>0$. Let $V^S(\Ga)$ be the set of stable vertices in~$V(\Ga)$. There are three types of unstable vertices:
\begin{gather*}
V^1(\Ga) = \{ v\in V(\Ga)\,|\, g_v=0, \,\val(v)=1, \,n_v=0\},\\
V^{1,1}(\Ga) = \{ v\in V(\Ga)\,|\, g_v=0,\, \val(v)=n_v=1\},\\
V^2(\Ga) = \{ v\in V(\Ga)\,|\, g_v=0, \val(v)=2,\, n_v=0\}.
\end{gather*}
Then $V(\Ga)$ is the disjoint union of $V^1(\Ga)$, $V^{1,1}(\Ga)$, $V^2(\Ga)$, and $V^S(\Ga)$.

The set of stable f\/lags is def\/ined to be
\begin{gather*}
F^S(\Gamma) = \big\{(e,v)\in F(\Ga)\,|\, v\in V^S(\Ga)\big\}.
\end{gather*}
\end{Definition}

Given a decorated graph $\vGa=\big(\Ga,\vf,\vd,\vg,\vs\big)$, the curves $C_e$ and the maps $f|_{C_e}\colon C_e\to \ell_{\ep_e}\subset X$ are determined by $\vGa$. If $v\notin V^S(\Ga)$ then $C_v$ is a point. If $v\in V^S(\Ga)$ then $C_v$ is a curve, and $y(e,v):= C_e\cap C_v$ is a node of $C$ for $e\in E_v$,
\begin{gather*}
\bigl(C_v, \{ y(e,v)\colon e\in E_v\} \cup \{ x_i\,|\, i\in S_v\} \bigr)
\end{gather*}
is a $(\val(v)+n_v)$-pointed, genus $g_v$ curve, which represents a point in $\Mbar_{g_v,\val(v)+n_v}$. We call this moduli space $\Mbar_{g_v,E_v\cup S_v}$ instead of $\Mbar_{g_v, \val(v)+n_v}$ because we would like to label the marked points on $C_v$ by $E_v\cup S_v$ instead of $\{ 1,2, \ldots, \val(v)+n_v\}$. Then
\begin{gather*}
\cM_{\vGa}=\prod_{v\in V^S(\Ga)}\Mbar_{g_v, E_v\cup S_v}.
\end{gather*}
The automorphism group $A_{\vGa}$ for any point $[f\colon (C,x_1,\ldots,x_n)\to X]\in \cF_{\vGa}$ f\/its in the following short exact sequence of groups:
\begin{gather*}
1\to \prod_{e\in E(\Ga)} \bZ_{d_e} \to A_{\vGa}\to \Aut\big(\vGa\big)\to 1,
\end{gather*}
where $\bZ_{d_e}$ is the automorphism group of the degree $d_e$ morphism
\begin{gather*}
f|_{C_e}\colon \ C_e\cong \bP^1\to \ell_{\ep_e}\cong \bP^1,
\end{gather*}
and $\Aut(\vGa)$ is the automorphism group of the decorated graph $\vGa=\big(\Ga,\vf,\vd,\vg,\vs\big)$. There is a morphism $i_{\vGa}\colon \cM_{\vGa}\to \MgX$ whose image is the f\/ixed locus $\cF_{\vGa}$ associated to $\vGa\in \GgX$. The morphism $i_{\vGa}$ induces an isomorphism $[\cM_{\vGa}/A_{\vGa}]\cong \cF_{\vGa}$.

\subsection{Virtual tangent and normal bundles}
Given a decorated graph $\vGa\in \GgX$ and a stable map $f\colon (C,x_1,\ldots,x_n)\to X$ which represents a point in the f\/ixed locus~$\cF_{\vGa}$ associated to~$\vGa$, let
\begin{gather*}
B_1 = \Hom(\Omega_C(x_1+\cdots+x_n),\cO_C),\qquad B_2 = H^0(C,f^*TX),\\
B_4 = \Ext^1(\Omega_C(x_1+\cdots+ x_n),\cO_C),\qquad B_5= H^1(C,f^*TX).
\end{gather*}
Then $B_1$, $B_2$, $B_4$, $B_5$ are representations of the torus $T$, so there is a direct sum decomposition $B_i = B_i^f \oplus B_i^m$, where $B_i^f \subset B_i$ is the $T$-invariant subspace. We have the following exact sequences:
\begin{gather*}
0\to B_1^f\to B_2^f\to T^{1,f}\to B_4^f\to B_5^f\to T^{2,f}\to 0, \\
 0\to B_1^m\to B_2^m\to T^{1,m}\to B_4^m\to B_5^m\to T^{2,m}\to 0.
\end{gather*}

The irreducible components of $C$ are
\begin{gather*}
\big\{ C_v\,|\, v\in V^S(\Ga)\big\}\cup\{C_e\,|\, e\in E(\Ga) \}.
\end{gather*}
The nodes of $C$ are
\begin{gather*}
\big\{ y_v=C_v \,|\, v\in V^2(\Ga) \big\} \cup \big\{ y(e,v) \,|\, (e,v)\in F^S(\Ga)\big\}.
\end{gather*}

\subsubsection{Automorphisms of the domain} \label{sec:aut}
Given any $(e,v)\in F(\Ga)$, let $y(e,v)=C_e\cap C_v$, and
def\/ine
\begin{gather*}
w_{(e,v)}:=e^T(T_{y(e,v)}C_e)=\frac{\bw(\ep_e,\si_v)}{d_e} \in H_T^2(y(e,v);\bQ)= M\otimes_\bZ \bQ.
\end{gather*}
We have
\begin{gather*}
B_1^f = \bigoplus_\edge \Hom(\Omega_{C_e}(y(e,v)+y(e,v')),\cO_{C_e})\\
\hphantom{B_1^f}{} = \bigoplus_\edge H^0(C_e, TC_e(-y(e,v)-y(e,v')),\\
B_1^m = \bigoplus_\vone T_{y(e,v)}C_e.
\end{gather*}

\subsubsection{Deformations of the domain} \label{sec:deform}
Given any $v\in V^S(\Ga)$, def\/ine a divisor $\bx_v$ of $C_v$ by
\begin{gather*}
\bx_v=\sum_{i\in S_v} x_i + \sum_{e\in E_v} y(e,v).
\end{gather*}
Let $\bL_{(e,v)}$ be the line bundle over $\Mbar_{g_v, E_v\cup S_v}$ whose f\/iber at the moduli point $[C_v,\bx_v]$ is the cotangent line $T^*_{y(e,v)}C_v$. Let
\begin{gather*}
\psi_{(e,v)}= c_1(\bL_{(e,v)}) \in A^1\big(\Mbar_{g_v,E_v\cup S_v}\big).
\end{gather*}
The torus $T$ acts trivially on $\Mbar_{g_v, E_v\cup S_v}$ and $\bL_{(e,v)}$, so $\psi_{(e,v)}$ can also be viewed as the $T$-equivariant f\/irst Chern class of $\bL_{(e,v)}$.

We have
\begin{gather*}
B_4^f= \bigoplus_{v\in V^S(\Ga)} \Ext^1(\Omega_{C_v}(\bx_v),\cO_C) = \bigoplus_{v\in V^S(\Ga)} T_{(C_v,{\bf x}_v)}\Mbar_{g_v, E_v\cup S_v},\\
B_4^m = \bigoplus_{v\in V^2(\Ga), E_v=\{e,e'\} } T_{y_v}C_e\otimes T_{y_v} C_{e'} \oplus \bigoplus_{(e,v)\in F^S(\Ga)} T_{y(e,v)}C_v\otimes T_{y(e,v)} C_e,
\end{gather*}
where
\begin{gather*}
e^T \big(T_{y_v}C_e \otimes T_{y_v} C_{e'}\big)= w_{(e,v)}+w_{(e',v)},\qquad v\in V^2(\Ga),\\
e^T \big(T_{y(e,v)}C_v \otimes T_{y(e,v)} C_e\big) = w_{(e,v)}-\psi_{(e,v)},\qquad v\in V^S(\Ga).
\end{gather*}

\subsubsection{Unifying stable and unstable vertices}
From the discussion in Sections \ref{sec:aut} and~\ref{sec:deform},
\begin{gather}
\frac{e^T(B_1^m)}{e^T(B_4^m)}=\prod_\vone w_{(e,v)} \prod_{v\in V^2(\Ga), \, E_v=\{e,e'\} }\frac{1}{w_{(e,v)}+ w_{(e',v)} }\nonumber\\
\hphantom{\frac{e^T(B_1^m)}{e^T(B_4^m)}=}{}\times \prod_{v\in V^S(\Ga)}\frac{1}{\prod\limits_{e\in E_v}(w_{(e,v)}-\psi_{(e,v)}) }. \label{eqn:Bonefour}
\end{gather}

Recall that
\begin{gather*}
\cM_\vGa =\prod_{v\in V^S(\Ga)} \Mbar_{g_v, E_v\cup S_v}.
\end{gather*}

To unify the stable and unstable vertices, we use the following convention for the empty sets $\Mbar_{0,1}$ and $\Mbar_{0,2}$. Let $w_1$, $w_2$ be formal variables.
\begin{enumerate}\itemsep=0pt
\item[(i)] $\Mbar_{0,1}$ is a $-2$-dimensional space, and
\begin{gather}\label{eqn:one}
\int_{\Mbar_{0,1}}\frac{1}{w_1-\psi_1}=w_1.
\end{gather}
\item[(ii)] $\Mbar_{0,2}$ is a $-1$-dimensional space, and
\begin{gather}\label{eqn:two}
\int_{\Mbar_{0,2}}\frac{1}{(w_1-\psi_1)(w_2-\psi_2)}= \frac{1}{w_1+w_2},\\
\label{eqn:one-one}
\int_{\Mbar_{0,2}}\frac{1}{w_1-\psi_1} =1.
\end{gather}
\item[(iii)] $\cM_{\vGa}=\prod\limits_{v\in V(\Ga)} \Mbar_{g_v, E_v\cup S_v}$.
\end{enumerate}

With the above conventions (i), (ii), (iii), we may rewrite \eqref{eqn:Bonefour} as
\begin{gather*}
\frac{e^T(B_1^m)}{e^T(B_4^m)}= \prod_{v\in V(\Ga)}\frac{1}{\prod\limits_{e\in E_v}(w_{(e,v)}-\psi_{(e,v)}) }.
\end{gather*}

The following lemma shows that the conventions (i) and (ii) are consistent with the stable case $\Mbar_{0,n}$, $n\geq 3$.
\begin{Lemma}\label{lemma:psi}
For any positive integer $n$ and formal variables $w_1,\ldots,w_n$, we have
\begin{gather*}
(a) \quad \int_{\Mbar_{0,n}}\frac{1}{\prod_{i=1}^n(w_i-\psi_i)} =\frac{1}{w_1\cdots w_n}(\frac{1}{w_1}+\cdots \frac{1}{w_n})^{n-3},\\
(b) \quad \int_{\Mbar_{0,n}}\frac{1}{w_1-\psi_1} = w_1^{2-n}.
\end{gather*}
\end{Lemma}
\begin{proof} (a) The cases $n=1$ and $n=2$ follow from the def\/initions~\eqref{eqn:one} and~\eqref{eqn:two}, respectively. For $n\geq 3$, we have
\begin{gather*}
\int_{\Mbar_{0,n}}\frac{1}{\prod\limits_{i=1}^n (w_i-\psi_i)} =\frac{1}{w_1\cdots w_n}\int_{\Mbar_{0,n}} \frac{1}{\prod\limits_{i=1}^n\big(1-\frac{\psi_i}{w_i}\big)}\\
\hphantom{\int_{\Mbar_{0,n}}\frac{1}{\prod\limits_{i=1}^n (w_i-\psi_i)}}{} =\frac{1}{w_1\cdots w_n}\sum_{a_1+\cdots+ a_n=n-3} w_1^{-a_1}\cdots w_n^{-a_n} \int_{\Mbar_{0,n}}\psi_1^{a_1}\cdots \psi_n^{a_n},
\end{gather*}
where
\begin{gather*}
\int_{\Mbar_{0,n}}\psi_1^{a_1}\cdots \psi_n^{a_n} =\frac{(n-3)!}{a_1!\cdots a_n!}.
\end{gather*}
So
\begin{gather*}
 \int_{\Mbar_{0,n}}\frac{1}{\prod\limits_{i=1}^n (w_i-\psi_i)}=\frac{1}{w_1\cdots w_n} \left(\frac{1}{w_1}+\cdots \frac{1}{w_n}\right)^{n-3}.
\end{gather*}

(b) The cases $n=1$ and $n=2$ follow from the def\/initions \eqref{eqn:one} and \eqref{eqn:one-one}, respectively. For $n\geq 3$, we have
\begin{gather*}
\int_{\Mbar_{0,n}}\frac{1}{w_1-\psi_1}
=\frac{1}{w_1}\int_{\Mbar_{0,n}} \frac{1}{1-\frac{\psi_1}{w_1}} =\frac{1}{w_1} w_1^{3-n} = w_1^{2-n}.\tag*{\qed}
\end{gather*}\renewcommand{\qed}{}
\end{proof}

\subsubsection{Deformation of the map}
Consider the normalization sequence
\begin{gather*}%\label{eqn:normalize}
0 \to \cO_C\to \bigoplus_{v\in V^S(\Ga)} \cO_{C_v} \oplus \bigoplus_{e\in E(\Ga)} \cO_{C_e} \to \bigoplus_{v\in V^2(\Ga)} \cO_{y_v} \oplus \bigoplus_{(e,v)\in F^S(\Ga)} \cO_{ y(e,v) }\to 0.
\end{gather*}
We twist the above short exact sequence of sheaves by~$f^*TX$. The resulting short exact sequence gives rise a long exact sequence of cohomology groups
\begin{gather*}
0 \to B_2 \to \bigoplus_{v\in V^S(\Ga)} H^0(C_v)\oplus \bigoplus_{e\in E(\Ga)}H^0(C_e) \to \bigoplus_{v\in V^2(\Ga)} T_{f(y_v)}X
\oplus \bigoplus_{(e,v)\in F^S(\Ga)} T_{f(y(e,v))}X \\
\hphantom{0}{} \to B_5 \to \bigoplus_{v\in V^S(\Ga)} H^1(C_v)\oplus \bigoplus_{e\in E(\Ga)}H^1(C_e) \to 0,
\end{gather*}
where
\begin{gather*}
H^i(C_v) = H^i\big(C_v, (f|_{C_v})^*TX\big) \cong H^i\big(C_v, \cO_{C_v}\big)\otimes T_{p_{\si_v}}X,\\
H^i(C_e) = H^i\big(C_e, (f|_{C_e})^*TX\big)
\end{gather*}
for $i=0,1$. We have
\begin{gather*}
H^0(C_v)= T_{p_{\si_v}}X,\\
H^1(C_v)= H^0\big(C_v,\omega_{C_v}\big)^\vee\otimes T_{p_{\si_v}}X.
\end{gather*}

\begin{Lemma} Let $\si\in V(\Up)$, so that $p_\si$ is a $T$-fixed point in $X$. Define
\begin{gather*}
\bw(\si)= e^T(T_{p_\si}X) \in H^{2r}_T(\pt;\bQ),\\
\bh(\si,g)= \frac{e^T(\bE^\vee\otimes T_{p_\si}X) }{e^T(T_{p_\si}X)}\in H^{2r(g-1)}_T\big(\Mbar_{g,n};\bQ\big).
\end{gather*}
Then
\begin{gather}\label{eqn:wsi}
\bw(\si)=\prod_{\ep\in E_\si} \bw(\ep,\si), \\ \label{eqn:hsig}
\bh(\si,g)=\prod_{\ep\in E_\si} \frac{\Lambda_g^\vee(\bw(\ep,\si))}{\bw(\ep,\si)},
\end{gather}
where $\Lambda^\vee_g(u)=\sum\limits_{i=0}^g (-1)^i \lambda_i u^{g-i}$.
\end{Lemma}
\begin{proof} $T_{p_\si}X =\bigoplus\limits_{\ep\in E_\si} T_{p_\si} \ell_\ep$, where $e^T(T_{p_\si}\ell_\ep)= \bw(\ep,\si)$. So
\begin{gather*}
e^T(T_{p_\si}X) = \prod_{\ep\in E_\si} \bw(\ep,\si),\\
\frac{e^T(\bE^\vee\otimes T_{p_\si}\ell_\ep)}{e^T(T_{p_\si}\ell_\ep)}
 = \prod_{\ep\in E_\si}\frac{e^T(\bE^\vee\otimes T_{p_\si}\ell_\ep)}{\bw(\ep,\si)},
\end{gather*}
where
\begin{gather*}
e^T\big(\bE^\vee\otimes T_{p_\si}\ell_\ep\big)= \sum_{i=0}^g (-1)^i c_i(\bE)c_1^T\big(T_{p_\si}\ell_\ep\big)^{g-i}= \sum_{i=0}^g (-1)^i\lambda_i \bw(\ep,\si)^{g-i}.\tag*{\qed}
\end{gather*}\renewcommand{\qed}{}
\end{proof}

The map $B_1\to B_2$ sends
\begin{gather*}
H^0(C_e, TC_e(-y(e,v)-y(e',v)))
\end{gather*}
 isomorphically to
 \begin{gather*}
 H^0\big(C_e, (f|_{C_e})^*T\ell_{\ep_e}\big)^f,
\end{gather*} the f\/ixed part of $H^0(C_e, (f|_{C_e})^*T\ell_{\ep_e})$.

\begin{Lemma}
Given $d\in \bZ_{>0}$ and $\ep\in E(\Up)_c$, define $\si$, $\si'$, $\ep_i$, $\ep_i'$, $a_i$ as in Section~{\rm \ref{GKM-graph}}, and let $f_d\colon \bP^1\to \ell_\ep\cong \bP^1$ be the unique degree $d$ map totally ramified over the two $T$-fixed points~$p_\si$ and~$p_{\si'}$ in~$\ell_\ep$. Define
\begin{gather*}
\bh(\ep,d)=\frac{e^T\big(H^1\big(\bP^1, f_d^* TX \big)^m \big)} {e^T \big(H^0\big(\bP^1, f_d^*TX\big)^m\big) }.
\end{gather*}
Then
\begin{gather}\label{eqn:etaud}
\bh(\ep,d)=\frac{(-1)^d d^{2d}}{ (d!)^2 \bw(\ep,\si)^{2d}}
\prod_{i=1}^{r-1} b\left(\frac{\bw(\ep,\si)}{d}, \bw(\ep_i,\si), da_i\right),
\end{gather}
where
\begin{gather}\label{eqn:b}
b(u,w,a)=\begin{cases}
\displaystyle \prod_{j=0}^a(w-ju)^{-1}, & a\in \bZ, \ a\geq 0,\\
\displaystyle \prod_{j=1}^{-a-1} (w+ju), & a\in \bZ, \ a<0.
\end{cases}
\end{gather}
\end{Lemma}

\begin{proof}
We use the notation in Section \ref{GKM-graph}. We have
\begin{gather*}
N_{\ell_\ep/X}=L_1\oplus \cdots \oplus L_{r-1}.
\end{gather*}
The weights of $T$-actions on $(L_i)_{p_\si}$ and $(L_i)_{p_\si}$ are $\bw(\ep_i,\si)$ and $\bw(\ep_i,\si)-a_i\bw(\ep,\si)$, respectively. The weights of $T$-actions on $T_0\bP^1$, $T_\infty \bP^1$, $(f_d^*L_i)_0$, $(f_d^*L_i)_\infty$ are $u:=\frac{\bw(\ep,\si)}{d}$, $-u$, $w_i:= \bw(\ep_i,\si)$, $w_i- da_i u$, respectively. By \cite[Example~19]{Liu},
\begin{gather*}
\ch_T\big(H^0\big(\bP^1, f_d^*L_i\big)-H^1\big(\bP^1, f_d^*L_i\big)\big)
=\begin{cases}
\displaystyle \sum_{j=0}^{da_i} e^{w_i-ju}, & a_i\geq 0,\\
\displaystyle \sum_{j=1}^{-da_i-1}e^{w_i+ju}, & a_i<0.
\end{cases}
\end{gather*}

Note that $w_i+ju$ is nonzero for any $j\in \bZ$ since $w_i$ and $u$ are linearly independent for $i=1,\ldots,n-1$. So
\begin{gather*}
\frac{e^T\big(H^1\big(\bP^1,f_d^*L_i\big)\big)}{e^T\big(H^0\big(\bP^1,f_d^*L_i\big)\big)} =\frac{e^T\big(H^1\big(\bP^1,f_d^*L_i\big)^m\big)}{e^T\big(H^0\big(\bP^1,f_d^*L_i\big)^m\big)} =b(u, w_i, da_i),
\end{gather*}
where $b(u,w,a)$ is def\/ined by \eqref{eqn:b}. By \cite[Example 19]{Liu},
\begin{gather*}
\ch_T\big(H^0\big(\bP^1, f_d^*T\ell_\ep\big)-H^1\big(\bP^1, f_d^*T\ell_\ep\big)\big)
=\sum_{j=0}^{2d} e^{du-ju}= 1+\sum_{j=1}^d \big(e^{j\bw(\ep,\si)/d}+ e^{-j\bw(\ep,\si)/d}\big).
\end{gather*}
So
\begin{gather*}
\frac{e^T\big(H^1\big(\bP^1,f_d^*T\ell_\ep\big)^m\big)}{e^T\big(H^0\big(\bP^1,f_d^*T\ell_\ep\big)^m\big)}
=\prod_{j=1}^d \frac{-d^2}{j^2 \bw(\ep,\si)^2 }=\frac{(-1)^d d^{2d}}{(d!)^2 \bw(\ep,\si)^{2d}}.
\end{gather*}
Therefore,
\begin{gather*}
 \frac{e^T\big(H^1\big(\bP^1, f_d^* TX \big)^m \big)} {e^T \big(H^0\big(\bP^1, f_d^*TX\big)^m\big) }
=\frac{e^T\big(H^1\big(\bP^1,f_d^*T\ell_\ep\big)^m\big)}{e^T\big(H^0\big(\bP^1,f_d^*T\ell_\ep\big)^m\big)}
\prod_{i=1}^{r-1}\frac{e^T\big(H^1\big(\bP^1,f_d^*L_i\big)^m\big)}{e^T\big(H^0\big(\bP^1,f_d^*L_i\big)^m\big)}\\
\hphantom{\frac{e^T\big(H^1\big(\bP^1, f_d^* TX \big)^m \big)} {e^T \big(H^0\big(\bP^1, f_d^*TX\big)^m\big) }}{}= \frac{(-1)^d d^{2d}}{(d!)^2 \bw(\ep,\si)^{2d}}\prod_{i=1}^{r-1}b\left(\frac{\bw(\ep,\si)}{d}, \bw(\ep_i,\si), da_i\right).\tag*{\qed}
\end{gather*}\renewcommand{\qed}{}
\end{proof}
Finally, $f(y_v)=p_{\si_v}= f(y(e,v))$, and
\begin{gather*}
e^T\big(T_{p_{\si_v}}X\big) =\bw(\si_v).
\end{gather*}
From the above discussion, we conclude that
\begin{gather*}
\frac{e^T(B_5^m)}{e^T(B_2^m)} = \prod_{v\in V^2(\Ga)} \bw(\si_v) \prod_{(e,v)\in F^S(\Ga)}\bw(\si_v)
 \prod_{v\in V^S(\Gamma)}\bh(\si_v,g_v) \prod_{e\in E(\Ga)}\bh(\ep_e,d_e)\\
\hphantom{\frac{e^T(B_5^m)}{e^T(B_2^m)}}{} = \prod_{v\in V(\Ga)}\bigl( \bh(\si_v,g_v) \bw(\si_v)^{\val(v)} \bigr)
 \prod_{e\in E(\Ga)}\bh(\ep_e,d_e),
\end{gather*}
where $\bw(\si)$, $\bh(\si,g)$, and $\bh(\ep,d)$ are def\/ined by
\eqref{eqn:wsi}, \eqref{eqn:hsig}, \eqref{eqn:etaud}, respectively.

\subsection{Contribution from each graph} \label{sec:each-graph}
\subsubsection{Virtual tangent bundle} We have $B_1^f=B_2^f$, $B_5^f=0$. So
\begin{gather*}
T^{1,f}=B_4^f =\bigoplus_{v\in V^S(\Ga)}T_{(C_v,{\bf x}_v)}\Mbar_{g_v, E_v\cup S_v},\qquad T^{2,f}=0.
\end{gather*}
We conclude that
\begin{gather*}
\bigg[\prod_{v\in V^S(\Ga)}\Mbar_{g_v, E_v\cup S_v}\bigg]^\vir =\prod_{v\in V^S(\Ga)}\big[\,\Mbar_{g_v, E_v\cup S_v} \big].
\end{gather*}
\subsubsection{Virtual normal bundle}
Let $N^\vir_{\vGa}$ be the pull back of the virtual normal bundle of $\cF_{\vGa}$ in $\MgX$ under $i_\vGa\colon \cM_{\vGa}\to \cF_{\vGa}$. Then
\begin{gather*}
\frac{1}{e^T\big(N^\vir_\vGa\big)} = \frac{e^T(B_1^m)e^T(B_5^m)}{e^T(B_2^m)e^T(B_4^m)}
=\prod_{v\in V(\Ga)}\frac{\bh(\si_v,g_v) \bw(\si_v)^{\val(v)} }{\prod\limits_{e\in E_v}(w_{(e,v)}-\psi_{(e,v)})}
 \prod_{e\in E(\Ga)}\bh(\ep_e,d_e).
\end{gather*}

\subsubsection{Integrand} Given $\si\in V(\Up)$, let
\begin{gather*}
i_\si^*\colon \ A^*_T(X)\to A^*_T(p_\si) =\bQ[u_1,\ldots,u_r]
\end{gather*}
be induced by the inclusion $i_\si\colon p_\si \to X$. Then
\begin{gather}
i_\vGa^*\prod_{i=1}^n\big(\ev_i^*\gamma_i^T \cup \big(\psi_i^T\big)^{a_i}\big)\nonumber\\
\qquad{} =\prod_{{v\in V^{1,1}(E)\atop S_v=\{i\}, \, E_v=\{e\}}} i^*_{\si_v}\gamma_i^T \big({-}w_{(e,v)}\big)^{a_i}
 \prod_{v\in V^S(\Ga)}\bigg(\prod_{i\in S_v} i^*_{\sigma_v}\gamma_i^T\prod_{e\in E_v}\psi_{(e,v)}^{a_i}\bigg).\label{eqn:integrand}
\end{gather}
To unify the stable vertices in $V^S(\Ga)$ and the unstable vertices in $V^{1,1}(\Ga)$, we use the following convention: for $a\in \bZ_{\geq 0}$,
\begin{gather}\label{eqn:one-one-a}
\int_{\Mbar_{0,2}}\frac{\psi_2^a}{w_1-\psi_1}=(-w_1)^a.
\end{gather}
In particular, \eqref{eqn:one-one} is obtained by setting $a=0$. With the convention \eqref{eqn:one-one-a}, we may re\-wri\-te~\eqref{eqn:integrand} as
\begin{gather*}
i_\vGa^*\prod_{i=1}^n\big(\ev_i^*\gamma_i^T \cup \big(\psi_i^T\big)^{a_i}\big)=
\prod_{v\in V(\Ga)}\bigg(\prod_{i\in S_v} i^*_{\sigma_v}\gamma_i^T\prod_{e\in E_v}\psi_{(e,v)}^{a_i}\bigg).
\end{gather*}

The following lemma shows that the convention \eqref{eqn:one-one-a} is consistent with
the stable case $\Mbar_{0,n}$, $n\geq 3$.
\begin{Lemma} Let $n$, $a$ be integers, $n\geq 2$, $a\geq 0$. Then
\begin{gather*}
\int_{\Mbar_{0,n}}\frac{\psi_2^a}{w_1-\psi_1}=\begin{cases}
\displaystyle \frac{\prod\limits_{i=0}^{a-1}(n-3-i)}{a!}w_1^{a+2-n}, & n=2\text{ or } 0\leq a\leq n-3,\\
0, &\text{otherwise}.
\end{cases}
\end{gather*}
\end{Lemma}
\begin{proof} The case $n=2$ follows from \eqref{eqn:one-one-a}. For $n\geq 3$,
\begin{gather*}
\int_{\Mbar_{0,n}}\frac{\psi_2^a}{w_1-\psi_1}=\frac{1}{w_1} \int_{\Mbar_{0,n}}\frac{\psi_2^a}{1-\frac{\psi_1}{w_1}}=w_1^{a+2-n} \int_{\Mbar_{0,n}} \psi_1^{n-3-a}\psi_2^a\\
\hphantom{\int_{\Mbar_{0,n}}\frac{\psi_2^a}{w_1-\psi_1}}{} = w_1^{a+2-n}\frac{(n-3)!}{(n-3-a)! a_!}=\frac{\prod\limits_{i=0}^{a-1}(n-3-i)}{a!} w_1^{a+2-n}.\tag*{\qed}
\end{gather*}\renewcommand{\qed}{}
\end{proof}

\subsubsection{Integral}\label{sec:vGa-integral}
The contribution of
\begin{gather*}
\int_{[\MgX^T]^{\vir,T}} \frac{i^*\prod\limits_{i=1}^n \big(\ev_i^*\gamma_i^T\cup \big(\psi_i^T\big)^{a_i}\big)}{e^T(N^\vir)}
\end{gather*}
from the f\/ixed locus $\cF_\vGa$ is given by
\begin{gather*}
\frac{1}{|A_{\vGa}|}\prod_{e\in E(\Ga)}\bh(\ep_e,d_e) \prod_{v\in V(\Ga)}\bigg(\bw(\si_v)^{\val(v)}\prod_{i\in S_v}i_{\si_v}^*\gamma_i^T\bigg) \prod_{v\in V(\Ga)}\int_{\Mbar_{g_v,E_v\cup S_v}}
\frac{\bh(\si_v,g_v) \prod\limits_{e\in E_v}\psi^{a_i}_{(e,v)} }{\prod\limits_{e\in E_v}(w_{(e,v)}-\psi_{(e,v)})},
\end{gather*}
where $|A_{\vGa}| =|\Aut(\vGa)| \prod\limits_{e\in E(\Ga)}d_e$.

\subsection{Sum over graphs}
Summing over the contribution from each graph $\vGa$ given in Section \ref{sec:vGa-integral} above, we obtain the following formula.
\begin{Theorem}\label{main}
\begin{gather}
\big \langle \tau_{a_1}\big(\gamma_1^T\big)\cdots \tau_{a_n}\big(\gamma_n^T\big)\big\rangle^{X_T}_{g,\beta} \nonumber\\
\qquad {} = \sum_{\vGa\in G_{g,n}(X,\beta)} \frac{1}{|\Aut(\vGa)|}
\prod_{e\in E(\Ga)} \frac{\bh(\ep_e,d_e)}{d_e}
\prod_{v\in V(\Ga)}\bigg(\bw(\si_v)^{\val(v)} \prod_{i\in S_v} i_{\si_v}^* \gamma_i^T \bigg) \nonumber\\
\qquad\quad {} \times \prod_{v\in V(\Ga)} \int_{\Mbar_{g,E_v \cup S_v} }
\frac{\bh(\si_v,g_v)\prod\limits_{i\in S_v} \psi_i^{a_i}}{\prod\limits_{e \in E_v} (w_{(e,v)}-\psi_{(e,v)})},\label{eqn:sum}
\end{gather}
where $\bh(\ep,d)$, $\bw(\si)$, $\bh(\si,g)$ are given by \eqref{eqn:etaud}, \eqref{eqn:wsi}, \eqref{eqn:hsig}, respectively, and we have the following convention for the $v\notin V^S(\Ga)$:
\begin{gather*}
\int_{\Mbar_{0,1}}\frac{1}{w_1-\psi_2}= w_1,\qquad \int_{\Mbar_{0,2}}\frac{1}{(w_1-\psi_1)(w_2-\psi_2)}=\frac{1}{w_1+w_2},\\
\int_{\Mbar_{0,2}}\frac{\psi_2^a}{w_1-\psi_1}=(-w_1)^a,\qquad a\in \bZ_{\geq 0}.
\end{gather*}
\end{Theorem}

Given $g\in \bZ_{\geq 0}$, $r$ weights $\vec{w}=\{ w_1,\ldots,w_r\}$, $r$ partitions $\vec{\mu} =\{ \mu^1,\ldots, \mu^r\}$, and $a_1,\ldots,a_k\in \bZ$,
let $\ell(\mu^i)$ be the length of $\mu^i$, and let $\ell(\vmu)=\sum\limits_{i=1}^r \ell(\mu^i)$. We def\/ine
\begin{gather*}
\langle \tau_{a_1},\ldots,\tau_{a_k}\rangle_{g,\vmu,\vec{w}} =\int_{\Mbar_{g,\ell(\vmu)+k}}
\prod_{i=1}^r\Bigg(\frac{\Lambda_g^\vee(w_i) w_i^{\ell(\vmu)-1} }{\prod\limits_{j=1}^{\ell(\mu^i)}\big(\frac{w_i}{\mu^i_j}-\psi^i_j\big)}\Bigg)
\prod_{b=1}^k\psi_b^{a_i}.
\end{gather*}
Given $v\in V(\Ga)$, def\/ine $\vec{w}(v)= \{ \bw(\ep,\si_v)\,|\, (\ep,\si_v)\in F(\Up)\}$. Given $v\in V(\Ga)$, and $\ep\in E_{\si_v}$, let $\mu^{v,\ep}$ be a~(possibly empty) partition def\/ined by $\{ d_e\,|\, e\in E_v,\, \vf(e)=\ep\}$, and def\/ine $\vmu(v)=\{ \mu^{v,\ep}\,|\, (\ep,\si_v)\in F(\Up)\}$. Then~\eqref{eqn:sum} can be rewritten as
\begin{gather*}
\big\langle \tau_{a_1}\big(\gamma_1^T\big)\cdots \tau_{a_n}\big(\gamma_n^T\big)\big\rangle^{X_T}_{g,\beta}\\
\qquad {} = \sum_{\vGa\in G_{g,n}(X,\beta)} \frac{1}{ |\Aut(\vGa)| }
\prod_{e\in E(\Ga)} \frac{\bh(\ep_e,d_e)}{d_e}
\prod_{v\in V(\Ga)}\bigg(\prod_{i\in S_v} i_{\si_v}^* \gamma_i
\bigg\langle \prod_{i\in S_v}\tau_{a_i}\bigg\rangle_{g_v,\vmu(v),\vec{w}(v)}\bigg).
\end{gather*}

Recall that
\begin{gather*}
g= \sum_{v\in V(\Ga)}g_v + |E(\Ga)| - |V(\Ga)| +1,
\end{gather*}
so
\begin{gather*}
2g-2 =\sum_{v\in V(\Ga)} (2g_v-2 +\val(v)).
\end{gather*}

Given $\vGa=\big(\Ga,\vf,\vd,\vg,\vs\big)$, let $\vGa'=\big(\Ga,\vf,\vd,\vs\big)$ be the decorated graph obtained by forgetting the genus map. Let $G_n(X,\beta)=\big\{ \vGa'\,|\, \vGa\in \cup_{g\geq 0}\GgX \big\}$. Def\/ine
\begin{gather}\label{eqn:gene}
\big\langle\tau_{a_1}\big(\gamma_1^T),\dots, \tau_{a_n}\big(\gamma_n^T\big)\,|\, u\big\rangle^{X_T}_\beta
= \sum_{g\geq 0} u^{2g-2}\big\langle\tau_{a_1}\big(\gamma_1^T\big),\dots, \tau_{a_n}\big(\gamma_n^T\big)\big\rangle^{X_T}_{g,\beta},
\\
\langle \tau_{a_1},\ldots,\tau_{a_k}\,|\, u \rangle_{\vmu,\vec{w}}=
\sum_{g\geq 0}u^{2g-2+\ell(\vmu)}\langle \tau_{a_1},\ldots,\tau_{a_k}\rangle_{g,\vmu,\vec{w}}.\nonumber
\end{gather}
Then we have the following formula for the generating function \eqref{eqn:gene}.
\begin{Theorem}\label{generating}
\begin{gather*}
\big\langle \tau_{a_1}\big(\gamma_1^T\big)\cdots \tau_{a_n}\big(\gamma_n^T\big)\,|\, u\big\rangle^{X_T}_\beta
= \sum_{\vGa'\in G_n(X,\beta)} \frac{1}{|\Aut(\vGa)|}
\prod_{e\in E(\Ga)} \frac{\bh(\ep_e,d_e)}{d_e} \\
\hphantom{\big\langle \tau_{a_1}\big(\gamma_1^T\big)\cdots \tau_{a_n}\big(\gamma_n^T\big)\,|\, u\big\rangle^{X_T}_\beta=}{} \times \prod_{v\in V(\Ga)}\bigg(\prod_{i\in S_v} i_{\si_v}^* \gamma_i^T
\bigg\langle \prod_{i\in S_v}\tau_{a_i}\,|\, u \bigg\rangle_{\vmu(v),\vec{w}(v)}\bigg).
\end{gather*}
\end{Theorem}

\subsection*{Acknowledgements}

The f\/irst author would like to thank Tom Graber for his suggestion of generalizing the computations for toric manifolds in~\cite{Liu} to GKM manifolds. The second author would like to thank the Columbia University for hospitality during his visits. We wish to thank Rahul Pandharipande for his comments on an earlier version of this paper. This work is partially supported by NSF DMS-1159416 and NSF DMS-1206667.

\pdfbookmark[1]{References}{ref}
\LastPageEnding


\begin{thebibliography}{99}
\footnotesize\itemsep=0pt

\bibitem{AtBo}
Atiyah M.F., Bott R., The moment map and equivariant cohomology,
 \href{https://doi.org/10.1016/0040-9383(84)90021-1}{\textit{Topology}} \textbf{23} (1984), 1--28.

\bibitem{Be1}
Behrend K., Gromov--{W}itten invariants in algebraic geometry, \href{https://doi.org/10.1007/s002220050132}{\textit{Invent.
Math.}} \textbf{127} (1997), 601--617, \mbox{\href{https://arxiv.org/abs/alg-geom/9601011}{alg-geom/9601011}}.

\bibitem{Be2}
Behrend K., Localization and {G}romov--{W}itten invariants, in Quantum
 Cohomology ({C}etraro, 1997), \href{https://doi.org/10.1007/978-3-540-45617-9_2}{\textit{Lecture Notes in Math.}}, Vol. 1776,
 Springer, Berlin, 2002, 3--38.

\bibitem{BeFa}
Behrend K., Fantechi B., The intrinsic normal cone, \href{https://doi.org/10.1007/s002220050136}{\textit{Invent. Math.}}
 \textbf{128} (1997), 45--88, \href{https://arxiv.org/abs/alg-geom/9601010}{alg-geom/9601010}.

\bibitem{BeMa}
Behrend K., Manin Yu., Stacks of stable maps and {G}romov--{W}itten invariants,
 \href{https://doi.org/10.1215/S0012-7094-96-08501-4}{\textit{Duke Math.~J.}} \textbf{85} (1996), 1--60, \href{https://arxiv.org/abs/alg-geom/9506023}{alg-geom/9506023}.

\bibitem{CLL}
Chen L., Li Y., Liu K., Localization, {H}urwitz numbers and the {W}itten
 conjecture, \href{https://doi.org/10.4310/AJM.2008.v12.n4.a5}{\textit{Asian~J. Math.}} \textbf{12} (2008), 511--518,
 \href{https://arxiv.org/abs/math.AG/0609263}{math.AG/0609263}.

\bibitem{DeMu}
Deligne P., Mumford D., The irreducibility of the space of curves of given
 genus, \href{https://doi.org/10.1007/BF02684599}{\textit{Inst. Hautes \'Etudes Sci. Publ. Math.}} \textbf{36} (1969), 75--109.

\bibitem{Fa}
Faber C., Algorithms for computing intersection numbers on moduli spaces of
 curves, with an application to the class of the locus of {J}acobians, in New
 Trends in Algebraic Geometry ({W}arwick, 1996), \href{https://doi.org/10.1017/CBO9780511721540.006}{\textit{London Math. Soc.
 Lecture Note Ser.}}, Vol.~264, Cambridge University Press, Cambridge, 1999,
 93--109, \href{https://arxiv.org/abs/alg-geom/9706006}{alg-geom/9706006}.

\bibitem{GKM}
Goresky M., Kottwitz R., MacPherson R., Equivariant cohomology, {K}oszul
 duality, and the localization theorem, \href{https://doi.org/10.1007/s002220050197}{\textit{Invent. Math.}} \textbf{131}
 (1998), 25--83.

\bibitem{GrPa}
Graber T., Pandharipande R., Localization of virtual classes, \href{https://doi.org/10.1007/s002220050293}{\textit{Invent.
 Math.}} \textbf{135} (1999), 487--518, \mbox{\href{https://arxiv.org/abs/alg-geom/9708001}{alg-geom/9708001}}.

\bibitem{GZ}
Guillemin V., Zara C., Equivariant de {R}ham theory and graphs,
 \href{https://doi.org/10.4310/AJM.1999.v3.n1.a3}{\textit{Asian~J. Math.}} \textbf{3} (1999), 49--76, \href{https://arxiv.org/abs/math.AG/9808135}{math.AG/9808135}.

\bibitem{Ha}
Hartshorne R., Algebraic geometry, \href{https://doi.org/10.1007/978-1-4757-3849-0}{\textit{Graduate Texts in Mathematics}},
 Vol.~52, Springer-Verlag, New York~-- Heidelberg, 1977.

\bibitem{HKKPTVVZ}
Hori K., Katz S., Klemm A., Pandharipande R., Thomas R., Vafa C., Vakil R.,
 Zaslow~E., Mirror symmetry, \textit{Clay Mathematics Monographs}, Vol.~1,
 Amer. Math. Soc., Providence, RI, Clay Mathematics Institute, Cambridge, MA,
 2003.

\bibitem{Ka}
Kazarian M., K{P} hierarchy for {H}odge integrals, \href{https://doi.org/10.1016/j.aim.2008.10.017}{\textit{Adv. Math.}}
 \textbf{221} (2009), 1--21, \href{https://arxiv.org/abs/0809.3263}{arXiv:0809.3263}.

\bibitem{KaL}
Kazarian M.E., Lando S.K., An algebro-geometric proof of {W}itten's conjecture,
 \href{https://doi.org/10.1090/S0894-0347-07-00566-8}{\textit{J.~Amer. Math. Soc.}} \textbf{20} (2007), 1079--1089,
 \href{https://arxiv.org/abs/math.AG/0601760}{math.AG/0601760}.

\bibitem{KiL}
Kim Y.-S., Liu K., Virasoro constraints and {H}urwitz numbers through asymptotic
 analysis, \href{https://doi.org/10.2140/pjm.2009.241.275}{\textit{Pacific~J. Math.}} \textbf{241} (2009), 275--284.

\bibitem{KnMu}
Knudsen F.F., Mumford D., The projectivity of the moduli space of stable
 curves. {I}.~{P}reliminaries on ``det'' and ``{D}iv'', \href{https://doi.org/10.7146/math.scand.a-11642}{\textit{Math. Scand.}}
 \textbf{39} (1976), 19--55.

\bibitem{Kn2}
Knudsen F.F., The projectivity of the moduli space of stable curves.
 {II}.~{T}he stacks {$M_{g,n}$}, \href{https://doi.org/10.7146/math.scand.a-12001}{\textit{Math. Scand.}} \textbf{52} (1983),
 161--199.

\bibitem{Kn3}
Knudsen F.F., The projectivity of the moduli space of stable curves.
 {III}.~{T}he line bundles on {$M_{g,n}$}, and a~proof of the projectivity of
 {$\overline M_{g,n}$} in characteristic~{$0$}, \href{https://doi.org/10.7146/math.scand.a-12002}{\textit{Math. Scand.}}
 \textbf{52} (1983), 200--212.

\bibitem{Ko1}
Kontsevich M., Intersection theory on the moduli space of curves and the matrix
 {A}iry function, \href{https://doi.org/10.1007/BF02099526}{\textit{Comm. Math. Phys.}} \textbf{147} (1992), 1--23.

\bibitem{Ko2}
Kontsevich M., Enumeration of rational curves via torus actions, in The Moduli
 Space of Curves ({T}exel {I}sland, 1994), \href{https://doi.org/10.1007/978-1-4612-4264-2_12}{\textit{Progr. Math.}}, Vol.~129,
 Birkh\"auser Boston, Boston, MA, 1995, 335--368, \href{https://arxiv.org/abs/hep-th/9405035}{hep-th/9405035}.

\bibitem{LiTi1}
Li J., Tian G., Virtual moduli cycles and {G}romov--{W}itten invariants of
 algebraic varieties, \href{https://doi.org/10.1090/S0894-0347-98-00250-1}{\textit{J.~Amer. Math. Soc.}} \textbf{11} (1998),
 119--174, \href{https://arxiv.org/abs/alg-geom/9602007}{alg-geom/9602007}.

\bibitem{Liu}
Liu C.-C.M., Localization in {G}romov--{W}itten theory and orbifold
 {G}romov--{W}itten theory, in Handbook of Moduli, {V}ol.~{II}, \textit{Adv.
 Lect. Math. (ALM)}, Vol.~25, Int. Press, Somerville, MA, 2013, 353--425,
 \href{https://arxiv.org/abs/1107.4712}{arXiv:1107.4712}.

\bibitem{Mi}
Mirzakhani M., Weil--{P}etersson volumes and intersection theory on the moduli
 space of curves, \href{https://doi.org/10.1090/S0894-0347-06-00526-1}{\textit{J.~Amer. Math. Soc.}} \textbf{20} (2007), 1--23.

\bibitem{MuZ}
Mulase M., Zhang N., Polynomial recursion formula for linear {H}odge integrals,
 \href{https://doi.org/10.4310/CNTP.2010.v4.n2.a1}{\textit{Commun. Number Theory Phys.}} \textbf{4} (2010), 267--293,
 \href{https://arxiv.org/abs/0908.2267}{arXiv:0908.2267}.

\bibitem{Mu}
Mumford D., Towards an enumerative geometry of the moduli space of curves, in
 Arithmetic and Geometry, {V}ol.~{II}, \href{https://doi.org/10.1007/978-1-4757-9286-7_12}{\textit{Progr. Math.}}, Vol.~36,
 Birkh\"auser Boston, Boston, MA, 1983, 271--328.

\bibitem{OP1}
Okounkov A., Pandharipande R., Gromov--{W}itten theory, {H}urwitz numbers, and
 matrix models, in Algebraic Geometry~-- {S}eattle 2005, {P}art~1,
 \href{https://doi.org/10.1090/pspum/080.1/2483941}{\textit{Proc. Sympos. Pure Math.}}, Vol.~80, Amer. Math. Soc., Providence, RI,
 2009, 325--414, \href{https://arxiv.org/abs/math.AG/0101147}{math.AG/0101147}.

\bibitem{Sp}
Spielberg H., A formula for the {G}romov--{W}itten invariants of toric
 varieties, Ph.D.~Thesis, {U}niversit\'e Louis Pasteur (Strasbourg I),
 Strasbourg, 1999, \href{https://arxiv.org/abs/math.AG/0006156}{math.AG/0006156}.

\bibitem{Wi}
Witten E., Two-dimensional gravity and intersection theory on moduli space, in
 Surveys in Dif\/ferential Geometry ({C}ambridge, {MA}, 1990), Lehigh
 University, Bethlehem, PA, 1991, 243--310.

\end{thebibliography}
\end{document}